\documentclass[12pt]{amsart}
\usepackage{amsmath, amsthm, amscd, amsfonts}

\makeatletter \oddsidemargin.2375in \evensidemargin \oddsidemargin
\theoremstyle{definition}

\theoremstyle{remark}

\numberwithin{equation}{section}
\usepackage{amssymb}
\usepackage{graphicx}
\usepackage{latexsym}

\begin{document}
%----------------
\title[]
{}
\author[]
{}
\thanks{$^{1}$M.Fallahpour@kiau.ac.ir}
\thanks{$^{2}$Corresponding author: M-Khodabin@kiau.ac.ir}
\thanks{$^{3}$Maleknejad@iust.ac.ir}
\date{}
\maketitle
%--------------------------
%------------------
\begin{minipage}{15cm}
\begin{center}
{\textbf{\Large Approximation solution of two-dimensional linear stochastic Volterra integral equation by applying the Haar wavelet}}\\
\end{center}
\end{minipage}
\vspace{2mm}
\begin{center} \index{M. Fallahpour}\index{M.
Khodabin}\index{K. Maleknejad}
{\textbf{M. Fallahpour$^{1}$, M. Khodabin$^{2}$, K. Maleknejad$^{3}$}\\

{\small Department of Mathematics, Karaj Branch, Islamic Azad University, Karaj, Iran.\\

\vspace{0.2cm}} }
\end{center}
\vspace{5mm}

%--------------------------
\begin{center}
\begin{minipage}{11.5cm}
\textbf{Abstract.} {\small Numerical solution of one-dimensional
stochastic integral equations because of the randomness has its
own problems, i.e. some of them no have analytical solution or
finding their analytic solution is very difficult. This problem
for two-dimensional equations is twofold. Thus, finding an
efficient way to approximate solutions of these equations is an
essential requirement. To begin this important issue in this
paper, we will give an efficient method based on Haar wavelet to
approximate a solution for the two-dimensional linear stochastic
Volterra integral equation. We also give an example to demonstrate
the accuracy of the method.}
\end{minipage}
\end{center}\vspace{1cm}
\begin{center}
\begin{minipage}{15cm}
\textbf{Mathematical subject classification:} {\small 65C30, 65C20, 60H20, 60H35, 68U20.}\\
\textbf{\small Keywords:} {\small Haar wavelet; Two-dimensional
stochastic Volterra integral equation; Brownian motion process;
Ito integral.}
\end{minipage}
\end{center}
\vspace{1cm}
\section{\textbf{Introduction}}
As we know, two dimensional ordinary integral equations provide an
important tool for modeling a numerous problems in engineering and
science $ [6,7] $. The second kind of two-dimensional integral
equations may arise from some problems of
nonhomogeneous elasticity and electrostatics $ [8] $. \\

Dobner presented an equivalent formulation of the Dorboux problem
as a two-dimensional Volterra integral equation $ [9] $. We can
also see this kind of equations in contact problems for bodies
with complex properties $ [10,11] $, and in the theory of radio
wave propagation $ [12] $, and in the theory of the elastic
problem of axial translation of a rigid elliptical disc-inclusion
$ [13] $, and various physical, mechanical and biological
problems. Some numerical schemes have been
inspected for resolvent of two-dimensional ordinary integral
equations by several probers. Computational complexity of
mathematical operations is the most important obstacle for solving
ordinary integral equations in higher dimensionas.

The Nystrom method $ [14] $, collocation method $ [15,16,17] $,
Gauss product quadrature rule method $ [18] $, Galerkin method $ [19] $, using triangular
fuctions $ [20,21] $, Legender polynomial method $ [22] $,
differential transform method $ [23] $, meshless method $ [24] $,
Bernstein polynomials method $ [25] $ and Haar wavelet method $
[26] $. This paper is first focused on proposing a generic
framework for numerical solution of two-dimensional ordinary
linear Volterra integral equations of second kind. The use of the
Haar wavelet for the numerical solution of linear integral
equations has previously been discussed in $ [1] $ and references
therein. The paper $ [1] $ should be considered as a logical
continuation of the papers $ [2-4] $. In $ [2] $ a new numerical
method based on Haar wavelet is introduced for solution of
nonlinear one-dimensional Fredholm and Volterra integral
equations. In $ [3] $ the Haar wavelet method $ [2] $ is extended
to numerical solution of integro-differential equation. In $ [4] $
the Haar wavelet method $ [2,3] $ is improved in terms of
efficiency by introducing one-dimentional Haar wavelet
approximation of the kernel function.
The method $ [1] $ is fundamentally different from the other numerical methods based on Haar wavelet for the numerical solution of integral equations as it approximates kernel function using Haar wavelet.

The general hyperbolic differential equation is defined as $ [9] $
\begin{equation}
u_{xy}=r(x,y,u,u_{x},u_{y})\ \ , \ \ {(x,y)\in B}
\subseteq{\mathbb{R}^{2}}
\end{equation}
\[u_{1\hat{B}}=h(x,y) \ \ \ \ , \ \ \ \ \hat{B}\subseteq \partial B,\]
where the domain $ B $ and the subset $ \hat{B} $ of the border $ \partial B $ are chosen according to the different initial value problems.
It's easy to show that the integral form of $ (1.1) $ is given by two-dimentional Volterra integral equation
\[g(x,y)=f(x,y)+\int _{0}^{ y}\int _{0}^{x}K_{1}(x,y,s,t) g(s,t)dsdt.\]
Similarly, if we import statistical noise in to $ (1.1) $, we can obtain two-dimensional linear stochastic Volterra integral equation of the second kind, i.e.
\begin{equation}\label{GrindEQ__1_}
g(x,y)=f(x,y)+\int _{0}^{ y}\int _{0}^{x}K_{1}(x,y,s,t) g(s,t)dsdt
\end{equation}
\[\ \ \ \ \ \ \ \ \ \ \ \ \ \ \ \ \ \ \ \ \ \ \ \ \ \ \ \ \ \ \ +\int _{0}^{y}\int _{0}^{x}K_{2}(x,y,s,t) g(s,t) dB(s)dB(t)\]
\[\ \ \ \ \ \ \ \ \ \ \ \ \ \ \ \ \ \ \ \ \ \ \ \ \ \ \ \ \ \ \ \ \ \ \ \ \ \ \ \ \ \ \ \ \ \ \ \ \ \ \ \ \ \ \ \ \ \ \ \ \ \ (x,y)\in[0,1]\times [0,1] \ \ \ , \ \ \ {s}\leqslant {x<t}\leqslant {y.}\]\\
where the kernels $ K_{1}(x,y,s,t) $ and $ K_{2}(x,y,s,t) $ in $
(1,2) $ are known functions and $ f(x,y) $ is also a known
function whereas $ g(x,y) $ is unknown function and is called the
solution of two-dimensional stochastic integral equation. The
condition $s\leqslant x<t\leqslant y$ is necessary for
adaptability to the filtration $\{F_t; 0\leq t \leq1\}$ where
$F_t=\sigma\{B(s); 0\leq s \leq 1\}$.\\

\textbf{Lemma 1.}Put $\phi(t,s)=K(x,y,s,t)g(s,t)$. Let $\phi$ be a
function in ${\textit{L}}^{2}([0,1]^{2})$. Then there exists a
sequence ${\phi_{n}}$ of off-diagonal step functions such that

\[\lim_{n\rightarrow\infty}\int_{a}^{b}\int_{a}^{b}\mid\phi(t,s)-\phi_{n}(t,s)\mid^{2}dt ds = 0.\]\\

\textbf{Definition 1.} Let $ \phi \in \textit{L}^{2}([0,1]^{2}) $.
Then the double Wiener-It$ \hat{•} $o integral of $ {\phi} $ is
defined as
\[\int_{a}^{b} \int_{a}^{b}\phi(t,s) dB(t)dB(s)=\lim_{n\rightarrow \infty}\int_{a}^{b}\int_{a}^{b}\phi_{n}(t,s) dB(t)dB(s)\ \ \ \ \ \ in \ \ \  \textit{L}^{2}(\Omega).\]\\

\section{\textbf{Haar Wavelets}}
\noindent A wavelet family
$(\psi_j,_i\left(y\right))_{j\in{N},i\in{Z}}$ is an orthonormal
subfamily of the Hilbert space $L^{2}(R)$ with the property that
all function in the wavelet family are generated from a fixed
function $\psi$ called mother wavelet through dilations and
translations.The wavelet family satisfies the following relation
\[\psi_j,_i\left(y\right)=2^{j/2}\psi\left(2^{j}y-i\right).\]
\noindent For Haar wavelet family on the interval $[0,1)$ we have:
\[h_{1}(y)=\left\{ \begin{array}{cc} 1, & for $y$\in[0,1)\\0, & otherwise,\end{array}\right.\]\\
\[h_{i}(y)=\left\{\begin{array}{cc} 1 & for$y$ \in[\alpha,\beta)\ \ \ \ \ \ \ \ \ \\ -1,  &for$y$ \in[\beta,\gamma)\ \ \ \ \ \ \ \ \ \ \\0, &\ \ \ \ \ \ otherwise,\ \ \ \ \  i=2,3, ... , \end{array}\right.\]\\
\noindent\[\alpha=\frac{n}{m},\ \ \ \ \beta=\frac{(n+0.5)}{m},\ \ \ \ \gamma=\frac{(n+1)}{m};\]
\[\ m=2^{\ell},\ \ \ \ \ \ \ell=0, 1, ... ,\ \ \ \ \ \ n=0, 1, ... , m-1.\]\\
\noindent The integer $ \ell $ indicats the level of the wavelet
and $ n $ is the translation parameter. Any square integrable
function $ f(y) $ defined on $ [0,1) $ can be expressed as
follows:
\[f(y)=\sum _{i=1}^{\infty }a_{i} h_{i} (y),\]
where $ a_{i} $ are real constants.\\
For approximation aim we consider a maximum value $ L $ of the integer $ \ell, $ level of the Haar wavelet in the above definition. The integer $ L $ is then called maximum level of resolution. We also define integer $ M=2^{L}$. Hence for any square integrable function $ f(y) $ we have a finite sum of Haar wavelets as follows:
\[f(y)\thickapprox\sum _{i=1}^{2M }a_{i} h_{i} (y).\]
\noindent The following notation is introduced $ [1] $:
\begin{equation} \label{GrindEQ__1_}
p_{i,1}(y)=\int _{0}^{y}h_{i}(u)du,
\end{equation}where by the definition of Haar wavelet equation $ (2.1) $ reduce to
\[p_{i,1}(y)=\left\{\begin{array}{cc} y-\alpha, & \ \ \ \ \ \ \ for$ y $ \in[\alpha,\beta)\ \ \ \ \ \ \ \\ \gamma-y, & \ \ \ \ \ \ \ for$ y $\in[\beta,\gamma)\ \ \ \ \ \ \ \\ 0, & \ \ \ elsewhere.\ \ \ \ \ \ \ \end{array}\right.\]\\
We can also the following stochastic notation introduce
\begin{equation} \label{GrindEQ__1_}
q_{i,1}(y)=\int _{0}^{y}h_{i}(u)dB(u),
\end{equation}
where equation $ (2.2) $ can be evaluated similarly by the definition of Haar wavelet and is given as follow:\\

\[q_{i,1}(y)=\left\{\begin{array}{cc} B(y)-B(\alpha), & \ \ \ \ \ \ for $ y $ \in[\alpha,\beta)\ \ \ \ \ \ \ \ \ \ \ \ \ \ \\ 2B(\beta)-B(\alpha)-B(y), & \ \ \ \ \ \ \ for$ y $ \in[\beta,\gamma)\ \ \ \ \ \ \ \ \ \ \ \ \ \ \ \ \\ 0,& \ \ \ elsewhere.\ \ \ \ \ \ \ \ \ \ \ \ \ \ \end{array}\right.\]\\

\section{\textbf{Numerical method}}

\noindent In this section, proposed numerical method $[1]$ will be discussed for two-dimensional linear stochastic Volterra integral equation of the second kind. In the first subsection, we state some results for efficient evalution of two-dimensional Haar wavelet approximations. In the second subsection, we apply these results for finding numerical solutions equation $ (1.1) $.\\
For Haar wavelet approximation of a function $ f(x,y) $ of two
real variables $ x $ and $ y $, we assume that the domain $ 0\leq
x,y \leq1 $ is divided into a grid of size $ 2M\times2N $ using
the following collocation points
\begin{equation} \label{GrindEQ__1_}
x_{m}=\frac{m-0.5}{2M}, m=1,2, ... ,2M,
\end{equation}
\begin{equation} \label{GrindEQ__1_}
y_{n}=\frac{n-0.5}{2N}, n=1,2, ... ,2N.
\end{equation}\\

\textbf{3.1\ \ \ \ \ Two-dimensional Haar wavelet system}\\

\noindent A real-valued function $ G(x,y) $ of two real variables $ x $ and $ y $ can be approximated using two-dimensional Haar wavelets basis as $ [1,27] $:
\begin{equation} \label{GrindEQ__1_}
G(x,y)\thickapprox\sum _{p=1}^{2M}\sum _{q=1}^{2N}b_{p,q} h_{p} (x)h_{q} (y).
\end{equation}
In order to calculate the unknown coefficients {$ b_{i,j} $}'s,
the collocation points defined in Eqs. $ (3.1) $ and $ (3.2) $ are
substituted in Eq. $ (3.3) $.Hence, we obtain the following $
2M\times2N $ linear system with unknowns {$ b_{i,j} $}'s:
\begin{equation}
G(x_{m},y_{n})=\sum _{p=1}^{2M}\sum _{q=1}^{2N}b_{p,q} h_{p} (x_{m})h_{q} (y_{n}), m=1, 2, ... , 2M,\ \ \ n=1, 2, ... , 2N.
\end{equation}
The solution of system $ (3.4) $ can be calculated from the following theorem.\\
\textbf{Theorem 2.} The solution of the system $ (3.4) $ is given below:\\
\[b_{1,1}=\frac{1}{2M\times2N}\sum _{p=1}^{2M}\sum _{q=1}^{2N}G(x_{m},y_{n}),\]
\[b_{i,1}=\frac{1}{\rho_{1}\times2N}\left(\sum _{p=\alpha_{1}}^{\beta_{1}}\sum _{q=1}^{2N}G(x_{m},y_{n})-\sum _{p=\beta_{1}+1}^{\gamma_{1}}\sum _{q=1}^{2N}G(x_{m},y_{n})\right) \ \ \ , \ \ \ i=2, 3, ... , 2M,\]
\[b_{1,j}=\frac{1}{2M\times\rho_{2}}\left(\sum _{p=1}^{2M}\sum _{q=\alpha_{2}}^{\beta_{2}}G(x_{m},y_{n})-\sum _{p= 1}^{2M}\sum _{q=\beta_{2}+1}^{\gamma_{2}}G(x_{m},y_{n})\right)\ \ \ , \ \ \ j=2, 3, ... , 2N,\]
\[b_{i,j}=\frac{1}{\rho_{1}\times\rho_{2}}\left(\sum _{p=\alpha_{1}}^{\beta_{1}}\sum _{q=\alpha_{2}}^{\beta_{2}}G(x_{m},y_{n})-\sum _{p=\alpha_{1}}^{\beta_{1}}\sum _{q=\beta_{2}+1}^{\gamma_{2}}G(x_{m},y_{n})-\sum _{p=\beta_{1}+1}^{\gamma_{1}}\sum _{q=\alpha_{2}}^{\beta_{2}}G(x_{m},y_{n})\right.\]
\[\left.\ \ \ \ \ \ +\sum _{p=\beta_{1}+1}^{\gamma_{1}}\sum _{q=\beta_{2}+1}^{\gamma_{2}}G(x_{m},y_{n})\,\right)\ \ \ , \ \ \ i=2, 3, ... , 2M\ \ , \ \ j=2, 3, ... , 2N,\]
where\\
\begin{equation}
\begin{array}{cc} \alpha_{1}=\rho_{1}(\sigma_{1}-1)+1,&\\
\ \
\beta_{1}=\rho_{1}(\sigma_{1}-1)+\frac{\rho_{1}}{2},&\\\gamma_{1}=\rho_{1}\sigma_{1},\
\ \ \ \ \ \ \ \ \ \ \\\rho_{1}=\frac{2M}{\tau_{1},}\ \ \ \ \ \
\ \ \ \ \ \ \ \ \\\sigma_{1}=i-\tau_{1},\ \ \ \ \ \ \ \ \ \ \
&\\\tau_{1}=2^{\lfloor\log _{2}(i-1)\rfloor}\ \ \ \ \ \
&\end{array}.
\end{equation}
and similarly,
\begin{equation}
\begin{array}{cc} \alpha_{2}=\rho_{2}(\sigma_{2}-1)+1,\\ \ \ \beta_{2}=\rho_{2}(\sigma_{2}-1)+\frac{\rho_{2}}{2},\\\gamma_{2}=\rho_{2}\sigma_{2},\ \ \ \ \ \ \ \ \ \ \ \ \\\rho_{2}=\frac{2N}{\tau_{2}},\ \ \ \ \ \ \ \ \ \ \ \ \ \ \\\sigma_{2}=j-\tau_{2},\ \ \ \ \ \ \ \ \ \ \ \\\tau_{2}=2^{\lfloor\log _{2}(j-1)\rfloor}\ \ \ \ \ \ \end{array}.
\end{equation}
\textit{Proof. }See $ [2] $.\[\ \ \ \ \ \ \ \ \ \ \ \ \ \ \ \ \ \ \ \ \ \ \ \ \ \ \ \ \ \ \ \ \ \ \ \ \ \ \ \ \ \ \ \ \ \ \ \ \ \ \ \ \ \ \ \ \ \ \ \ \ \ \ \ \ \ \ \ \ \ \ \ \ \ \ \ \ \ \ \ \ \ \ \ \ \ \ \ \ \ \ \ \ \ \ \ \ \ \ \ \ \ \ \ \ \ \ \ \ \ \ \ \ \ \ \ \ \ \ \ \ \  \ \ \ \ \ \ \ \ \ \ \ \ \ \ \ \ \ \ \square\]\\

Consider a function $ G(x,y,s,t) $ of four variables $ x,y,s $ and $ t $. Suppose $ G(x,y,s,t) $ is approximated using two-dimensional Haar wavelet as follows $ [1] $:\\
\begin{equation} \label{GrindEQ__1_}
G(x,y,s,t)\thickapprox\sum _{p=1}^{2M}\sum _{q=1}^{2N}b_{p,q}(x,y) h_{p}(s)h_{q}(t).
\end{equation}\\
Substituting the collocation points
\[s_{i}=\frac{i-0.5}{2M}\ \ , \ \ i=1,2, ... ,2M,\]
and
\[t_{j}=\frac{j-0.5}{2N}\ \ , \ \ j=1,2, ... ,2N,\]
we obtain the linear system
\begin{equation}
G(x,y,s_{i},t_{j})\thickapprox\sum _{p=1}^{2M}\sum _{q=1}^{2N}b_{p,q}(x,y)h_{p} (s_{i})h_{q} (t_{j})\ \ , \ \ i=1,2, ... ,2M\ \ , \ \ j=1,2, ... ,2N.
\end{equation}\\
\textbf{Corollary 1.} The solution of the system $ (3.8) $ for any
value of $ x.y\in[0,1] $ is given as follows $ [1] $:
\[b_{1,1}(x,y)=\frac{1}{2M\times2N}\sum _{p=1}^{2M}\sum _{q=1}^{2N}G(x,y,s_{p},t_{q}),\]
\[b_{i,1}(x,y)=\frac{1}{\rho_{1}\times2N}\left(\sum _{p=\alpha_{1}}^{\beta_{1}}\sum _{q=1}^{2N}G(x,y,s_{p},t_{q})-\sum _{p=\beta_{1}+1}^{\gamma_{1}}\sum _{q=1}^{2N}G(x,y,s_{p},t_{q})\right), i=2, 3, ... , 2M,\]
\[b_{1,j}(x,y)=\frac{1}{2M\times\rho_{2}}\left(\sum _{p=1}^{2M}\sum _{q=\alpha_{2}}^{\beta_{2}}G(x,y,s_{p},t_{q})-\sum _{p= 1}^{2M}\sum _{q=\beta_{2}+1}^{\gamma_{2}}G(x,y,s_{p},t_{q})\right),j=2, 3, ... , 2N,\]
\[b_{i,j}(x,y)=\frac{1}{\rho_{1}\times\rho_{2}}\left(\sum _{p=\alpha_{1}}^{\beta_{1}}\sum _{q=\alpha_{2}}^{\beta_{2}}G(x,y,s_{p},t_{q})-\sum _{p=\alpha_{1}}^{\beta_{1}}\sum _{q=\beta_{2}+1}^{\gamma_{2}}G(x,y,s_{p},t_{q})\ \ \ \ \ \ \ \ \ \ \ \ \ \ \ \ \ \ \ \ \ \ \ \ \ \ \ \ \ \ \ \ \ \ \ \ \ \ \ \ \ \ \ \ \ \ \ \ \right.\]
\[\left.\ \ \ \ \ \ \ \ \ \ \ \ \ \ \ \ \ \ -\sum _{p=\beta_{1}+1}^{\gamma_{1}}\sum _{q=\alpha_{2}}^{\beta_{2}}G(x,y,s_{p},t_{q})
+\sum _{p=\beta_{1}+1}^{\gamma_{1}}\sum _{q=\beta_{2}+1}^{\gamma_{2}}G(x,y,s_{p},t_{q})\right),\]
\[\ \ \ \ \ \ \ \ \ \ \ \ \ \ \ \ \ \ \ \ \ \ \ \ \ \ \ \ \ \ \ \ \ \ \ \ \ \ \ \ \ \ \ \ \ \ \ \ \ \ \ \ \ \ \ \ \ \ \ \ \ \ \ \ \ \ \ \ \ \ \ \ \ \ \ \ \ \ \ \ \ \ \ \ \ \ \ \ \ \ \ \ \ \ \ \ \ \ \ \ \ \ \ \ \ \  i=2, 3, ... , 2M\ \ , \ \ j=2, 3, ... ,2N,\]
where $ \alpha_{1},\beta_{1},\gamma_{1} $ and $ \rho_{1} $ are defined as in Eq. $ (3.5) $ and $ \alpha_{2},\beta_{2},\gamma_{2} $ and $ \rho_{2} $ are defined as in Eq. $ (3.6) $.\\
\textbf{Corollary 2.} Suppose a function $ G(x,y) $ of two variables $ x $ and $ y $ is approximated using Haar wavelet approximation given in Eq. $ (3.3) $. Suppose further that $ G(x,y) $ is known at collocation points $ (x_{m},y_{m}) $, $ m=1, 2, ... ,2M, n=1, 2, ... ,2N. $ Then the approximate value of the function $ G(x,y) $ at any other point of the domain can be calculated as follows $ [1] $:\\
\[G(x,y)=\frac{1}{2M\times2N}\sum _{p=1}^{2M}\sum _{q=1}^{2N}G(x_{m},y_{m})h_{1}(x)h_{1}(y)\ \ \ \ \ \ \ \ \ \ \ \ \ \ \ \ \ \ \ \ \ \ \ \ \ \ \ \ \ \ \ \ \ \ \ \ \ \ \ \ \ \ \ \ \ \ \ \ \ \ \ \ \ \ \ \ \ \ \ \ \ \ \ \ \ \ \ \ \ \ \ \ \ \ \ \ \ \ \ \ \ \ \ \ \ \ \ \ \ \ \ \ \]
\[\ \ \ \ \ \ \ +\sum _{i=1}^{2M}\frac{1}{\rho_{1}\times2N}\left(\sum _{p=\alpha_{1}}^{\beta_{1}}\sum _{q=1}^{2N}G(x_{m},y_{m})-\sum_{p=\beta_{1}+1}^{\gamma_{1}}\sum _{q=1}^{2N}G(x_{m},y_{m})\right)h_{i}(x)h_{1}(y)\]
\[\ \ \ \ \ \ \ \ \ +\sum _{j=1}^{2N}\frac{1}{2M\times\rho_{2}}\left(\sum _{p=1}^{2M}\sum_{q=\alpha_{2}}^{\beta_{2}}G(x_{m},y_{m})-\sum _{p=1}^{2M}\sum_{q=\beta_{2}+1}^{\gamma_{2}}G(x_{m},y_{m})\right)h_{1}(x)h_{j}(y)\]
\[\ \ \ \ \ \ \ \ \ \ \ \ \ \ \ \ \ \ \ \ \ \ \ \ \ \ \ \ +\sum _{i=1}^{2M}\sum _{j=1}^{2N}\frac{1}{\rho_{1}\rho_{2}}\left(\sum _{p=\alpha_{1}}^{\beta_{1}}\sum_{q=\alpha_{2}}^{\beta_{2}}G(x_{m},y_{m})-\sum_{p=\alpha_{1}}^{\beta_{1}}\sum_{q=\beta_{2}+1}^{\gamma_{2}}F(x_{m},y_{m})\right.\]
\[\ \ \ \ \ \ \ \ \ \ \ \ \ \ \ \ \ \ \ \ \ \ \ \ \ \ \ \ \ \ \ \ \ \ \left.-\sum _{p=\beta_{1}+1}^{\gamma_{1}}\sum_{q=\alpha_{2}}^{\beta_{2}}G(x_{m},y_{m})
+\sum_{p=\beta_{1}+1}^{\gamma_{1}}\sum_{q=\beta_{2}+1}^{\gamma_{2}}G(x_{m},y_{m})\right)h_{i}(x)h_{j}(y),\]\\
where $ \alpha_{1},\beta_{1},\gamma_{1} $ and $ \rho_{1} $ are defined as in Eq. $ (3.5) $ and $ \alpha_{2},\beta_{2},\gamma_{2} $ and $ \rho_{2} $ are defined as in Eq. $ (3.6) $.\\
\textbf{3.2\ \ \ \ Two-dimensional linear stochastic Volterra integral equation}\\

\noindent Consider the two-dimensional linear stochastic Volterra integral equation $ (1.2) $. Assume that the function $ K(x,y,s,t) g(s,t) $ is approximated using two-dimensional Haar wavelet as follows:
\begin{equation}
K_{1}(x,y,s,t)g(s,t)\thickapprox\sum _{i=1}^{2M}\sum _{j=1}^{2N}b_{i,j}(x,y) h_{i}(s)h_{j}(t).
\end{equation}
\begin{equation}
K_{2}(x,y,s,t)g(s,t)\thickapprox\sum _{i=1}^{2M}\sum _{j=1}^{2N}c_{i,j}(x,y) h_{i}(s)h_{j}(t).
\end{equation}
With this approximation Eq. $ (1.2) $ can be writen as follows:
\begin{equation}
g(x,y)=f(x,y)+\int _{0}^{ y}\int _{0}^{x}\sum _{i=1}^{2M}\sum _{j=1}^{2N}b_{i,j}(x,y) h_{i}(s)h_{j}(t)dsdt\ \ \ \ \ \ \ \ \ \ \ \ \ \ \ \ \ \ \ \ \ \ \ \ \ \ \ \
\end{equation}
\[\ \ \ \ \ \ \ \ \ \ \ \ \ \ \ \ \ \ \ +\int _{0}^{y}\int _{0}^{x}\sum _{i=1}^{2M}\sum _{j=1}^{2N}c_{i,j}(x,y) h_{i}(s)h_{j}(t) dB(s)dB(t).\]
Eq. $ (3.11) $ can be written in a more compact form using the notations introduced in equations $ (2.1) $ and $ (2.2) $ and is given as follows:
\[g(x,y)=f(x,y)+\sum _{i=1}^{2M}\sum _{j=1}^{2N}b_{i,j}(x,y)p_{i,1}(x)p_{j,1}(y)+\sum _{i=1}^{2M}\sum _{j=1}^{2N}c_{i,j}(x,y)q_{i,1}(x)q_{j,1}(y).\]
Substituting the collocation points given in $ (3.1) $ and $ (3.2) $, we obtain the following system of equations:
\[g(x_{m},y_{n})=f(x_{m},y_{n})+\sum _{i=1}^{2M}\sum _{j=1}^{2N}b_{i,j}(x_{m},y_{n})p_{i,1}(x_{m})p_{j,1}(y_{n})\]
\[\ \ \ \ \ \ \ \ \ \ \ \ \ \ \ \ \ \ \ \ \ \ \ \ \ \ \ \ \ \ \ \ \ \ \ \ \ \ \ \ \ \ \ \ \ \ \ \ \ \ \ \ \ \ \ \ \ \ +\sum _{i=1}^{2M}\sum _{j=1}^{2N}c_{i,j}(x_{m},y_{n})q_{i,1}(x_{m})q_{j,1}(y_{n}).\]
Now $ b_{i,j}\ \ , \ \ i=1, 2, ... , 2M\ \ , \ \ j=1, 2, ... , 2N $ and similarly $ c_{i,j}\ \ , \ \ i=1, 2, ... , 2M\ \ , \ \ j=1, 2, ... , 2N $ can be replaced with their expressions given in Corollary 1 and the following system of equations is obtained:
\begin{equation}\label{GrindEQ__1_}
g(x_{m},y_{n})=f(x_{m},y_{n})+\frac{p_{1,1}(x_{m})p_{1,1}(y_{n})}{2M\times2N}\sum _{p=1}^{2M}\sum _{q=1}^{2N}K_{1}(x_{m},y_{n},s_{p},t_{q}) g(s_{p},t_{q})+
\end{equation}
\[\sum _{i=2}^{2M}\frac{p_{i,1}(x_{m})p_{1,1}(y_{n})}{\rho_{1}\times 2N}\left(\sum _{p=\alpha_{1}}^{\beta_{1}}\sum _{q=1}^{2N}K_{1}(x_{m},y_{n},s_{p},t_{q})g(s_{p},t_{q})-\ \ \ \ \ \ \ \right.\]
\[\left.\ \ \ \ \ \ \ \ \ \ \ \ \ \ \ \ \ \ \ \ \ \ \ \ \ \ \ \ \ \ \ \ \ \ \ \ \ \ \ \ \ \ \ \ \ \ \ \ \ \ \ \ \ \ \ \ \ \ \sum_{p=\beta_{1}+1}^{\gamma_{1}}\sum_{q=1}^{2N}K_{1}(x_{m},y_{n},s_{p},t_{q})\ \ g(s_{p},t_{q})\right)+\]\[\sum _{j=2}^{2N}\frac{p_{1,1}(x_{m})p_{j,1}(y_{n})}{2M\times\rho_{2}}\left(\sum _{p=1}^{2M}\sum _{p=\alpha_{2}}^{\beta_{2}}K_{1}(x_{m},y_{n},s_{p},t_{q})\ \ g(s_{p},t_{q})-\ \ \ \ \ \ \ \right.\]\[\left.\ \ \ \ \ \ \ \ \ \ \ \ \ \ \ \ \ \ \ \ \ \ \ \ \ \ \ \ \ \ \ \ \ \ \ \ \ \ \ \ \ \ \ \ \ \ \ \ \ \ \ \ \ \ \ \ \ \ \sum _{p=1}^{2M}\sum_{q=\beta_{2}+1}^{\gamma_{2}}K_{1}(x_{m},y_{n},s_{p},t_{q})\\g(s_{p},t_{q})\right)+\]\[\sum _{i=2}^{2M}\sum_{j=2}^{2N}\frac{p_{i,1}(x_{m})p_{j,1}(y_{n})}{\rho_{1}\times\rho_{2}}\left(\sum _{p=\alpha_{1}}^{\beta_{1}}\sum _{q=\alpha_{2}}^{\beta_{2}}K_{1}(x_{m},y_{n},s_{p},t_{q})\\g(s_{p},t_{q})-\ \ \ \ \ \ \right.\]\[\left.\ \ \ \ \ \ \ \ \ \ \ \ \sum _{p=\alpha_{1}}^{\beta_{1}}\sum_{q=\beta_{2}+1}^{\gamma_{2}}K_{1}(x_{m},y_{n},s_{p},t_{q})\ \ g(s_{p},t_{q})-\sum_{p=\beta_{1}+1}^{\gamma_{1}}\sum _{q=\alpha_{2}}^{\beta_{2}}K_{1}(x_{m},y_{n},s_{p},t_{q})\\g(s_{p},t_{q})+\right.\]\[\left.\ \ \ \ \ \ \ \ \ \ \ \ \ \ \ \ \ \ \ \ \ \ \ \ \ \ \ \ \ \ \ \ \ \ \ \ \ \ \ \ \ \ \ \ \ \ \ \ \ \ \ \ \ \ \ \ \ \ \sum_{p=\beta_{1}+1}^{\gamma_{1}}\sum_{q=\beta_{2}+1}^{\gamma_{2}}K_{1}(x_{m},y_{n},s_{p},t_{q})\\g(s_{p},t_{q})\right)+\]

\[\frac{q_{1,1}(x_{m})q_{1,1}(y_{n})}{2M\times2N}\sum _{p=1}^{2M}\sum _{q=1}^{2N}K_{2}(x_{m},y_{n},s_{p},t_{q}) g(s_{p},t_{q})+\ \ \ \ \ \ \ \ \ \ \ \ \ \ \ \ \ \ \ \ \ \]

\[\sum _{i=2}^{2M}\frac{q_{i,1}(x_{m})q_{1,1}(y_{n})}{\rho_{1}\times 2N}\left(\sum _{p=\alpha_{1}}^{\beta_{1}}\sum _{q=1}^{2N}K_{2}(x_{m},y_{n},s_{p},t_{q})g(s_{p},t_{q})-\right.\]
\[\left.\ \ \ \ \ \ \ \ \ \ \ \ \ \ \ \ \ \ \ \ \ \ \ \ \ \ \ \ \ \ \ \ \ \ \ \ \ \ \ \ \ \ \ \ \ \ \ \ \ \ \ \ \ \ \ \ \ \ \ \ \sum_{p=\beta_{1}+1}^{\gamma_{1}}\sum_{q=1}^{2N}K_{2}(x_{m},y_{n},s_{p},t_{q})\ \ g(s_{p},t_{q})\right)+\]\[\sum _{j=2}^{2N}\frac{q_{1,1}(x_{m})q_{j,1}(y_{n})}{2M\times\rho_{2}}\left(\sum _{p=1}^{2M}\sum _{p=\alpha_{2}}^{\beta_{2}}K_{2}(x_{m},y_{n},s_{p},t_{q})\,\,g(s_{p},t_{q})-\ \ \ \ \ \ \ \ \ \ \ \ \ \ \ \ \ \ \right.\]\[\left.\ \ \ \ \ \ \ \ \ \ \ \ \ \ \ \ \ \ \ \ \ \ \ \ \ \ \ \ \ \ \ \ \ \ \ \ \ \ \ \ \ \ \ \ \ \ \ \ \ \ \ \ \ \ \ \ \ \ \ \ \sum _{p=1}^{2M}\sum_{q=\beta_{2}+1}^{\gamma_{2}}K_{2}(x_{m},y_{n},s_{p},t_{q})\ \ g(s_{p},t_{q})\right)+\]\[\sum _{i=2}^{2M}\sum_{j=2}^{2N}\frac{q_{i,1}(x_{m})q_{j,1}(y_{n})}{\rho_{1}\times\rho_{2}}\left(\sum _{p=\alpha_{1}}^{\beta_{1}}\sum _{q=\alpha_{2}}^{\beta_{2}}K_{2}(x_{m},y_{n},s_{p},t_{q})\ \ g(s_{p},t_{q})-\ \ \ \ \ \ \ \ \ \ \ \ \ \right.\]\[\left.\ \ \ \ \ \ \ \ \sum _{p=\alpha_{1}}^{\beta_{1}}\sum_{q=\beta_{2}+1}^{\gamma_{2}}K_{2}(x_{m},y_{n},s_{p},t_{q})\ \ g(s_{p},t_{q})-\sum_{p=\beta_{1}+1}^{\gamma_{1}}\sum _{q=\alpha_{2}}^{\beta_{2}}K_{2}(x_{m},y_{n},s_{p},t_{q})\ \ g(s_{p},t_{q})+\right.\]\[\left.\ \ \ \ \ \ \ \ \ \ \ \ \ \ \ \ \ \ \ \ \ \ \ \ \ \ \ \ \ \ \ \ \ \ \ \ \ \ \ \ \ \ \ \ \ \ \ \ \ \ \ \ \ \ \ \ \sum_{p=\beta_{1}+1}^{\gamma_{1}}\sum_{q=\beta_{2}+1}^{\gamma_{2}}K_{2}(x_{m},y_{n},s_{p},t_{q})\,\,g(s_{p},t_{q})\right),\]
\[\ \ \ \ \ \ \ \ \ \ \ \ \ \ \ \ \ \ \ \ \ \ \ \ \ \ \ \ \ \ \ \ \ \ \ \ \ \ \ \ \ \ \ \ \ \ \ \ \ \ \ \ \ \ \ \ \ \ \ \ \ \ \ m=1, 2, ... , 2M\ \ , \ \ n=1, 2, ... , 2N.\]
Eq. $ (3.12) $ represents $ 2M\times2N $ system which can be solved using either prevalent methods for solving linear systems. The solution of this system gives values of $ g(x,y) $ at the collocation points. The values of $ g(x,y) $ at points other than collocation points can be calculated using Corollary 2.\\

\section{\textbf{Numerical Example}}

\noindent In this section, the numerical example is given to demonstrate the applicability and accuracy of our method.
Consider the following linear 2D stochastic Volterra integral equation of second kind:
\[u(x,y)=f(x,y)+\int _{0}^{ y}\int _{0}^{x}(x+y+t-s) u(s,t)dsdt+\int _{0}^{y}\int _{0}^{x}(x+y+t+s) u(s,t) dB(s)dB(t)\]
where
\[f(x,y)=x+y-\dfrac{1}{12}xy(x^{3}+4x^{2}y+4xy^{2}+y^{3}).\]
The solutions mean together confidence interval at the collocation points for the present method for $ 1000 $ iterative of system $ (3.12) $ is shown in Table $ 1 $. In Figs. $ 1-4$, three-dimensional graphs of the approximate solution for various values of level $ L $ are shown.\\
\\

\textbf{Table\,\,1:}The solutions mean together confidence interval for above example\\

\begin{tabular}{|c|c|c|ccc|c|ccc |}
\hline
$J$ & $M$ & $2M$ & & $(x,y)$ & & $\bar{u}(x,y)$ & & \small{Confidence}& \\
&&&&&&&$L$&\small{Interval}& $U$\\
\hline
% after \\: \hline or \cline{col1-col2} \cline{col3-col4} ...
0 & 1 & 2 & &$(0.25,0.75)$ &&1.9951 &1.9951&&1.9951\\
\hline
& & & &$(0.125,0.375)$ &&1.06717 &1.06715 &&1.06719\\
1&2 &4 & &$(0.375,0.875)$ &&2.43481 &2.43194 &&2.43769\\
& & & &$(0.625,0.875)$ &&2.72591 &2.72036 &&2.73147\\
\hline
& & & &$(0.0625,0.4375)$ &&1.0498 &1.0498 &&1.04981\\
2&4 &8 & &$(0.3125,0.6875)$ &&2.15292 &2.15266 &&2.15318\\
& & & &$(0.8125,0.9375)$ &&2.68057 &2.67707 &&2.68407\\
\hline
& & & &$(0.03125,0.71875)$ &&1.51492 &1.51488 &&1.51495\\
3&8 &16 & &$(0.40625,0.53125)$ &&2.25846 &2.25532 &&2.26159\\
& & & &$(0.78125,0.96875)$ &&2.67147 &2.67139 &&2.67156\\
\hline
& & & &$(0.015625,0.609375)$ &&1.26132 &1.26129 &&1.26135\\
4&16 &32 & &$(0.296875,0.796875)$ &&2.25999 &2.25891 &&2.26107\\
& & & &$(0.859375,0.984375)$& &2.58876 &2.58870 &&2.58881\\
\hline
\end{tabular}~\\
\newpage
\begin{figure}[h]
\centering{
\includegraphics[scale=0.5,width=\linewidth,bb= 0 0 581 240]{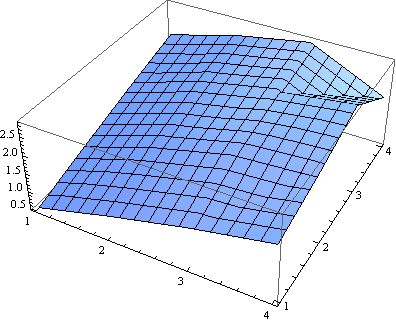} %\linewidth,scale=0.5]{Photo}%[width=10mm]{Photo.jpg}
\caption{Plot of approximate solution of level $ L=1 $ for test example. \label{overflow}}}
\end{figure}~\\~\\~\\~\\
\begin{figure}[h]
\centering{
\includegraphics[scale=0.5,width=\linewidth,bb= 0 0 581 240]{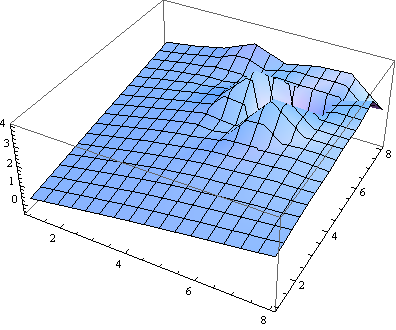} %\linewidth,scale=0.5]{Photo}%[width=10mm]{Photo.jpg}
\caption{Plot of approximate solution of level $ L=2 $ for test example. \label{overflow}}}
\end{figure}~\\~\\~\\~\\

\begin{figure}[h]
\centering{
\includegraphics[scale=0.5,width=\linewidth,bb= 0 0 581 240]{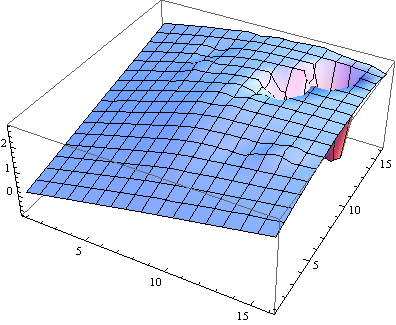} %\linewidth,scale=0.5]{Photo}%[width=10mm]{Photo.jpg}
\caption{Plot of approximate solution of level $ L=3 $ for test example. \label{overflow}}}
\end{figure}~\\~\\

\begin{figure}[h]
\centering{
\includegraphics[scale=0.5,width=\linewidth,bb= 0 0 581 240]{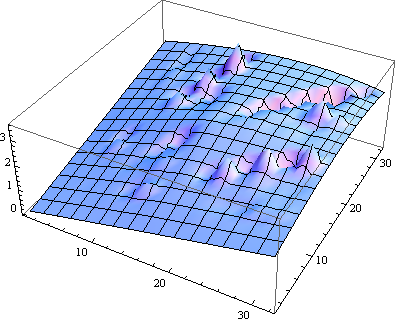} %\linewidth,scale=0.5]{Photo}%[width=10mm]{Photo.jpg}
\caption{Plot of approximate solution of level $ L=4 $ for test example. \label{overflow}}}
\end{figure}~\\ \\
\\

\newpage
\section{\textbf{Conclusion}}
As mentioned above, numerical solution of two-dimensional stochastic integral equations because of the randomness is very difficult or sometimes impossible. In this paper, we have successfully developed Haar wavelets numerical method for approximate a solution of two-dimensional linear stochastic Volterra integral equations. The example confirm that the method is considerably fast and highly accurate as sometimes lead to exact solution. Although, theoretically for getting higher accuracy we can set the method with larger values of M and N and also larger of the degree of approximation, p and q, but it leads to solving MN linear systems of size $ pq\times pq $, that have its difficulties. The method can be improved to be more accurate by using other numerical methods. Mathematica has been used for computations in this paper.\\

\end{document}